\documentclass[12pt]{article}

\def\sqr#1#2{{\vcenter{\vbox{\hrule height.#2pt
              \hbox{\vrule width.#2pt height#1pt \kern#1pt \vrule width.#2pt}
              \hrule height.#2pt}}}}
\def\signed #1{{\unskip\nobreak\hfil\penalty50
              \hskip2em\hbox{}\nobreak\hfil#1
              \parfillskip=0pt \finalhyphendemerits=0 \par}}
\def\endpf{\signed {$\sqr69$}}

\def\dbR{{\mathop{\rm l\negthinspace R}}}
\def\dbC{{\mathop{\rm l\negthinspace\negthinspace\negthinspace C}}}

\def\3n{\negthinspace \negthinspace \negthinspace }
\def\2n{\negthinspace \negthinspace }
\def\1n{\negthinspace }

\def\dbC{{\mathop{\rm l\negthinspace\negthinspace\negthinspace C}}}

\def\ds{\displaystyle}

\def\dbN{{\mathop{\rm l\negthinspace N}}}

\def\dbR{{\mathop{\rm l\negthinspace R}}}
\def\={\buildrel \triangle \over =}

\def\a{\alpha}
\def\b{\beta}

\def\d{\delta}
\def\e{\varepsilon}

\def\l{\lambda}

 \def\n{\nabla}

\def\t{\times}

\def\th{\theta}

\def\i{\infty}
\def\ns{\noalign{\ss} }

\def\G{\Gamma}

\def\Si{\Sigma}

\def\O{\Omega}

\def\cL{{\cal L}}

\def\cP{{\cal P}}

\def\cl{{\cal l}}
\def\no{\noindent}

\def\ss{\smallskip}
\def\ms{\medskip}

\def\q{\quad}
\def\qq{\qquad}
\def\hb{\hbox}

\def\pa{\partial}

\def\cd{\cdot}
\def\cds{\cdots}

\def\div{\hbox{\rm div$\,$}}

\def\cl{\overline}

\def\Re{{\mathop{\rm Re}\,}}
\def\Im{{\mathop{\rm Im}\,}}

\def\Sp{{\mathop{\rm Sp}\,}}
\def\|{\Big |}
\def\({\Big (}
\def\){\Big )}
\def\[{\Big[}
\def\]{\Big]}
\def\be{\begin{equation}}
\def\bel{\begin{equation}\label}
\def\ee{\end{equation}}
\def\bt{\begin{theorem}}
\def\bcd{\begin{condition}}
\def\ecd{\end{condition}}
\def\et{\end{theorem}}
\def\bc{\begin{corollary}}
\def\ec{\end{corollary}}
\def\bde{\begin{definition}}
\def\ede{\end{definition}}
\def\bl{\begin{lemma}}
\def\el{\end{lemma}}
\def\bp{\begin{proposition}}
\def\ep{\end{proposition}}
\def\br{\begin{remark}}
\def\er{\end{remark}}
\def\ba{\begin{array}}
\def\ea{\end{array}}
\def\ed{\end{document}}
\def\ns{\noalign{\ms}}
\def\ds{\displaystyle}

\def\square#1{\vbox{\hrule\hbox{\vrule height#1%
     \kern#1\vrule}\hrule}}
\def\rectangle#1#2{\vbox{\hrule\hbox{\vrule height#1%
     \kern#2\vrule}\hrule}}

\font\tenbb=msbm10 \font\sevenbb=msbm7 \font\fivebb=msbm5

\newfam\bbfam
\scriptscriptfont\bbfam=\fivebb \textfont\bbfam=\tenbb
\scriptfont\bbfam=\sevenbb

\def\oO{{\overline \O}}
\def\ov{{\overline v}}

\newtheorem{lemma}{Lemma}[section]
\newtheorem{remark}{Remark}[section]

\newtheorem{theorem}{Theorem}[section]
\newtheorem{corollary}{Corollary}[section]

\newtheorem{definition}{Definition}[section]
\newtheorem{proposition}{Proposition}[section]
\newtheorem{condition}{Condition}[section]

\makeatletter
   
   \@addtoreset{equation}{section}
\makeatother

\begin{document}
\title{\bf Logarithmic decay of hyperbolic equations with arbitrary boundary damping}

\author{Xiaoyu Fu\thanks{School of Mathematics,  Sichuan University,
Chengdu 610064,  China. This work was partially supported by the NSF
of China under grants 10525105 and 10771149. The author gratefully
acknowledges Dr. Kim-Dang Phung and Prof. Xu Zhang for their helps.
{\small\it E-mail:} {\small\tt rj\_xy@163.com}.}}

\date{}

\maketitle

\begin{abstract}

\no In this paper, we study the logarithmic stability for the
hyperbolic equations by arbitrary boundary observation. Based on
Carleman estimate, we first prove an estimate of the resolvent
operator of such equation. Then we prove the logarithmic stability
estimate for the hyperbolic equations without any assumption on an
observation subboundary.

\end{abstract}

\ms

\no{\bf 2000 Mathematics Subject Classification}.  Primary 93B05;
Secondary 93B07, 35B37.

\ms

\no{\bf Key Words}.  Logarithmic stability, hyperbolic equations,
Carleman estimate, resolvent operator.

\section{Introduction and main result}

Let $\O\subset\dbR^{n}$ be a bounded domain with boundary $\pa\O$ of
class $C^2$. Denote by $\nu=(\nu_1,\cdots,\nu_n)$ the unit outward
normal field along the boundary $\pa \O$, and $\oO$ the closure of
$\O$. For simplicity, in the sequel, we use the notation $ u_j={{\pa
u}\over{\pa x_j}}, $ where $x_j$ is the $j$-th coordinate of a
generic point $x=(x_1,\cdots,x_n)$ in $\dbR^n$. In a similar manner,
we use the notation $ w_j$, $v_j$, etc. for the partial derivatives
of $w$ and $v$ with respect to $x_j$. By $\overline c$ we denote the
complex conjugate of $c\in \dbC$. Throughout this paper, we will use
$C$ to denote a generic positive constant which may vary from line
to line (unless otherwise stated).

Let $a^{jk}(\cd)\in C^2(\cl{\O};\;\dbR)$ be fixed satisfying
 \bel{s3}
 a^{jk}=a^{jk}(x)=a^{kj}(x),\qq\forall\; x\in\cl \O,~j,k=1,2,\cdots,n,
 \ee
and for some constant $\b>0$,
 \bel{s4}
 \sum_{j,k=1}^na^{jk}(x)\xi^j\overline{\xi}^k
 \geq\b|\xi|^{2},\qq \forall\;(x,\xi)
 \in\overline{\O}\t \dbC^{n},
 \ee
where  $\xi=(\xi^1,\cdots,\xi^n)$. Define a formal differential
operator $\cP$ (associated with the matrix
$\big(a^{jk}(\cd)\big)_{n\times n}$) as follows:
  \bel{ca}
  \cP \=\sum_{j,k=1}^n\pa_k(a^{jk}\pa_j).
  \ee

Fix a real valued function $a(\cd)\in C^1(\pa\O;\;\dbR^+)$. In what
follows, we assume that
 \bel{a1}
\G_0\=\{x\in\pa\O;\;a(x)>0\}\neq\emptyset.
 \ee

The main purpose of this article is to study the logarithmic decay
of the following hyperbolic equations with a boundary damping term
$a(x)u_t$:
 \bel{0a1}\left\{\ba{ll}\ds
 u_{tt}-\cP u=0 &\hb{ in } \dbR^+\t \O,\\
 \ns\ds
\sum_{j,k=1}^na^{jk}u_j\nu_k=0&\hb{ on } \dbR^+\t \pa\O\setminus\G_0,\\
\ns\ds
\sum_{j,k=1}^na^{jk}u_j\nu_k+a(x)u_t=0&\hb{ on } \dbR^+\t \G_0,\\
\ns\ds (u(0), u_t(0))=(u^0,u^1)&\hb{ in } \O.
 \ea\right.\ee
Very interesting logarithmic decay results were given in \cite{C,
LR} for the above system under the regularity assumption that
$a^{jk}(\cd)$, $a(\cd)$ and $\pa\O$ are $C^\infty$-smooth (\cite{LR}
considered the special case $(a^{jk})_{n\t n}=I$, the identity
matrix). Note that, since the sub-boundary $\G_\d$ in which the
damping $a(x)u_t$ is (uniformly) effective may be very ``small", the
``geometric optics condition" introduced in \cite{BLR} is not
guaranteed for system (\ref{0a1}), and therefore, in general, one
can not expect exponential stability of this system. On the other
hand, as pointed in \cite{C, LR}, for some special case of system
(\ref{0a1}), logarithmic stability is the best decay rate.

Put $H\=H^1(\O)\t L^2(\O)$. Define an unbounded operator $A: \; H\to
H$ by (Recall that $ u_j^0={{\pa u^0}\over{\pa x_j}} $ )
 \bel{0a2}\left\{\ba{ll}
 A\=\left(\ba{cc}
 0&I\\\cP&0
 \ea
 \right),\\
 \ns\ds D(A)\=\left\{u=(u^0,u^1)\in H;\  Au\in H\right.,\\
 \ns\ds \qq\qq\qq\qq\left. \sum_{j,k=1}^na^{jk}u^0_{j}\nu_k\|_{ \pa\O\setminus\G_0}=0 ,\  \(\sum_{j,k=1}^na^{jk}u^0_{j}\nu_k+au^1\)\|_{\G_0}=0 \right\}. \ea\right.\ee
It is easy to show that $A$ generates a group
$\{e^{tA}\}_{t\in\dbR}$ on $H$.

The main result of this paper is stated as follows:
 \bt\label{0t1}
Let $a^{jk}(\cd)\in C^2(\cl{\O};\;\dbR)$ satisfy
(\ref{s3})--(\ref{s4}) and $a(\cd)\in C^1(\pa\O;\;\dbR^+)$ satisfy
(\ref{a1}). Then solutions $e^{tA}(u^0,u^1)\equiv (u,u_t)\in
C(\dbR;\; D(A))\cap C^1(\dbR;\;H)$ of system (\ref{0a1}) satisfy
  \bel{0a3}
 ||e^{tA}(u^0,u^1)||_{H}\le {C\over \ln (2+t)}||(u^0,u^1)||_{D(A)},\qq\forall\;(u^0,u^1)\in D(A),\ \forall\; t>0.
  \ee
 \et

Following \cite{B} (see also \cite{C}), Theorem \ref{0t1} is a
consequence of the following resolvent estimate for operator $A$:

 \bt\label{0t2} Under the assumptions in Theorem \ref{0t1}, there exists a constant $C>0$ such that

 i) if $\l\in
 \Sp(A)\setminus\{0\}$, then
  $$
 \Re\l<-{e^{-C|\Im\l|}\over C};
  $$

 ii) if
 $$
 \Re\l\in\left[-{e^{-C|\Im\l|}\over C},0\right],
 $$
 then
  $$
 ||(A-\l I)^{-1}||_{\cL(H)}\le Ce^{C|\Im\l|}, \q\hb{ for } |\l|>1.
  $$
 \et

We shall develop an approach based on global Carleman estimate to
prove Theorem \ref{0t2}, which is the main novelty of this paper.
Our approach, stimulated by \cite{LRS} (see also \cite{DZZ, Fu, Zh2,
Z1}), is different from that in \cite{B}, which instead employed the
classical local Carleman estimate and therefore needs
$C^\infty$-regularity for the data.

It would be quite interesting to establish better decay rate (than
logarithmic decay) for system (\ref{0a1}) under further conditions
(without geometric optics condition). There are some impressive
results in this respect, say \cite{BH, LR1, P1, P2} for polynomial
decay of system (\ref{0a1}) with special geometries. However, to the
best of the author's knowledge, the full picture of this problem is
still unclear. We refer to \cite{D, RZZ, ZZ2} for related works.

The rest of this paper is organized as follows. In section \ref{2},
we collect some useful preliminary results which will be useful
later. Another key preliminary, global Carleman estimate for
elliptic equations without inhomogeneous boundary condition, is
established in section \ref{ss5}. Sections \ref{ss6}--\ref{ss7} are
addressed to the proof of our main results.

\section{Some preliminaries}\label{2}

In this section, we collect some preliminaries which will be used in
the sequel.

To begin with, we recall the following result (which is an easy
consequence of known result in \cite{FI, WW}, for example).

\bl\label{0l4} There exists a function $\hat\psi\in C^2(\overline
\O)$ such that
 \bel{as}
 \left\{
 \ba{ll}\ds
 \hat\psi> 0 & \hb{ in } \O,\\
 \ns
 \ds |\n\hat\psi|>0& \hb{ in } \cl\O,\\
 \ns
 \ds
 \sum_{j,k=1}^na^{jk}\hat\psi_j\nu_k\le0  & \hb{ on }\pa\O\setminus\G_0.
 \ea
 \right.
 \ee
 \el

\ms

Next, for $n\in\dbN$, we denote by $O(\mu^n)$ a function of order
$\mu^n$ for large $\mu$ (which is independent of $\l$); by
$O_\mu(\l^n)$ a function of order $\l^n$ for fixed $\mu$ and for
large $\l$. We now show the following pointwise estimate, which is a
consequence of \cite[Theorem 2.1]{Fu} (see also \cite{Fu0}).
 \bl\label{l1} Let $a^{jk}\in
C^{2}(\dbR^{1+n};\;\dbR)$ satisfying (\ref{s3}). Assume $z\in
C^2(\dbR^{1+n};\;\dbC)$, $\Psi\in C^2(\dbR^{1+n};\dbR)$ and $\ell\in
C^4(\dbR^{1+n};\dbR)$. Set
 \bel{dw}
\th=e^\ell,\q v=\th z,\q
\Psi=-2\ell_{ss}-2\sum_{j,k=1}^n(a^{jk}\ell_j)_k.\
 \ee
 Then
 \bel{af1}\ba{ll}\ds
\th^2\|z_{ss}+\sum_{j,k=1}^n(a^{jk}z_j)_k\|^2+M_s+\div V\\
\ns\ds\ge
2\(3\ell_{ss}+\sum_{j,k=1}^n(a^{jk}\ell_j)_k\)|v_s|^2+4\sum_{j,k=1}^na^{jk}\ell_{js}(v_k\ov_s+\ov_kv_s)\\
\ns\ds\qq+\sum_{j,k=1}^nc^{jk}(v_k\ov_j+\ov_k v_j)+B|v|^2,
 \ea\ee
where
 \bel{asf}\left\{\ba{ll}\ds
 A=\ell_s^2+\sum_{j,k=1}^na^{jk}\ell_j\ell_k-\ell_{ss}-\sum_{j,k=1}^n(a^{jk}\ell_j)_k-\Psi,\\
 \ns\ds
 M=2\ell_s(|v_s|^2-\sum_{j,k=1}^na^{jk}\ov_jv_k)+2\sum_{j=1}^na^{jk}\ell_j(\ov_sv_j+v_s\ov_j)\\
 \ns\ds\qq\q-\Psi(\ov_sv+v_s\ov)+(2A\ell_s+\Psi_s)|v|^2,\\
 \ns\ds
 V=[V_1,\cdots,V_k,\cdots,V_n],\\
 \ns\ds
V_k=\sum_{j,j',k'=1}^n\Big\{-2a^{jk}\ell_j|v_s|^2+2a^{jk}\ell_s(\overline v_jv_s+v_j\overline v_s)-\Psi a^{jk}(v_j\ov+\ov_jv)\\
 \ns\ds\qq\q
+\(2a^{jk'}a^{j'k}-a^{jk}a^{j'k'}\)\ell_j(v_{j'}\ov_{k'}+\ov_{j'}v_{k'})+a^{jk}(2A\ell_j+\Psi_j-2a\ell_j\ell_t)|v|^2\Big\},\\
\ns\ds
c^{jk}=\sum_{j',k'=1}^n\[2(a^{j'k}\ell_{j'})_{k'}a^{jk'}-a^{jk}_{k'}a^{j'k'}\ell_{j'}+a^{jk}(a^{j'k'}\ell_{j'})_{k'}\]+a^{jk}\ell_{ss},\\
 \ns\ds
 B=\sum_{j,k=1}^n(a^{jk}\Psi_k)_j+2(A\ell_s)_s+2\sum_{j,k=1}^n(Aa^{jk}\ell_j)_k+2A\Psi.
 \ea\right.\ee
In particular,  for any function $\psi\in C^4(\dbR^{1+n};\;\dbR)$
satisfying $\psi_{sj}=0$ (j=1,\dots,n), and any $\l,\mu>1$, choosing
the function $\ell(s,x)$ to be
 \bel{c1}
\ell=\l\phi,\q \phi= e^{\mu\psi},
 \ee
then
  \bel{c3}
 \ba{ll}\ds
\hb{\rm Left hand side of (\ref{af1})}\ge
2\[\l\mu^2\phi\sum_{j,k=1}^na^{jk}\psi_j\psi_k+\l\phi
O(\mu)\]\(|v_s|^2+\sum_{j,k=1}^na^{jk}v_j\ov_k\)\\
\ns\ds\qq\qq\qq\qq\qq\q+2\[\l^3\mu^4\phi^3\|\sum_{j,k=1}^na^{jk}\psi_j\psi_k\|^2
+\l^3\phi^3O(\mu^3)+O_\mu(\l^2)\]|v|^2.
 \ea
 \ee
 \el

{\it Proof.} Using Theorem 2.1 in \cite{Fu} with $m=1+n$, and
 $$
 t=s,\q (a^{jk})_{m\times m}=\left(\ba{cc}
 1&0\\0&(a^{jk})_{n\times n}
 \ea
 \right).
 $$
By a direct calculation, we obtain (\ref{af1}).

On the other hand, by (\ref{c1}) and note that $\psi_{sj}=0$
($j=1,\dots,n$), it is easy to check that
 \bel{as6}\left\{\ba{ll}\ds
\ell_s=\l\mu\phi\psi_s,\qq\ell_j=\l\mu\phi\psi_j,\\
\ns\ds
\ell_{ss}=\l\mu^2\phi\psi_s^2+\l\mu\phi\psi_{ss},\q\ell_{jk}=\l\mu^2\phi\psi_j\psi_k+\l\mu\phi\psi_{jk},\q\ell_{js}=\l\mu^2\phi\psi_s\psi_j.
 \ea\right.\ee
Next, recalling the definition of $c^{jk}$ in (\ref{asf}), by
(\ref{as6}) and note that $a^{jk}$ satisfies (\ref{s3}), we have
 \bel{as10}\ba{ll}\ds
2\(3\ell_{ss}+\sum_{j,k=1}^n(a^{jk}\ell_j)_k\)|v_s|^2+4\sum_{j,k=1}^na^{jk}\ell_{js}(v_k\ov_s+\ov_kv_s)+\sum_{j,k=1}^nc^{jk}(v_k\ov_j+\ov_k
v_j)\\
\ns\ds=2\Big\{\l\mu^2\phi\[3|\psi_s|^2+\sum_{j,k=1}^na^{jk}\psi_j\psi_k\]+\l\phi
O(\mu)\Big\}|v_s|^2+8\l\mu^2\phi\sum_{j,k=1}^na^{jk}\psi_j\psi_sv_k\ov_s\\
\ns\ds\q+4\l\mu^2\|\sum_{j,k=1}^na^{jk}\psi_j\ov_k\|^2+2\Big\{\l\mu^2\phi\[\sum_{j,k=1}^na^{jk}\psi_j\psi_k+|\psi_s|^2\]+\l\phi
O(\mu)\Big\}\sum_{j,k=1}^na^{jk}v_k\ov_j\\
\ns\ds=4\l\mu^2\phi\|\psi_sv_s+\sum_{j,k=1}^na^{jk}\psi_jv_k\|^2+4\l\mu^2\|\sum_{j,k=1}^na^{jk}\psi_j\ov_k\|^2\\
\ns\ds\q+2\Big\{\l\mu^2\phi\[\sum_{j,k=1}^na^{jk}\psi_j\psi_k+|\psi_s|^2\]+\l\phi
O(\mu)\Big\}\(|v_s|^2+\sum_{j,k=1}^na^{jk}v_j\ov_k\)\\
\ns\ds\ge 2\[\l\mu^2\phi\sum_{j,k=1}^na^{jk}\psi_j\psi_k+\l\phi
O(\mu)\]\(|v_s|^2+\sum_{j,k=1}^na^{jk}v_j\ov_k\).
 \ea\ee
Further, by (\ref{as6}) and recalling (\ref{asf}) and (\ref{dw}) for
the definition of $A$ and $\Psi$, respectively, we have
 \bel{fd}\left\{\ba{ll}\ds
 \Psi=2\l\mu^2\phi\[|\psi_s|^2+\sum_{j,k=1}^na^{jk}\psi_j\psi_k\]+\l\phi O(\mu),\\
 \ns\ds A=(\l^2\mu^2\phi^2+\l\mu^2\phi)\[|\psi_s|^2+\sum_{j,k=1}^na^{jk}\psi_j\psi_k\]+\l\phi O(\mu).
 \ea\right.\ee
Therefore, by (\ref{asf}), and note that $a^{jk}$ satisfies
(\ref{s4}), we have
 \bel{as11}\ba{ll}\ds
B&\ds=2\l^3\mu^4\phi^3\|\sum_{j,k=1}^na^{jk}\psi_j\psi_k+|\psi_s|^2\|^2+\l^3\phi^3O(\mu^3)+O_\mu(\l^2)\\
\ns&\ds\ge2\l^3\mu^4\phi^3\|\sum_{j,k=1}^na^{jk}\psi_j\psi_k\|^2
+\l^3\phi^3O(\mu^3)+O_\mu(\l^2).
 \ea\ee
Combining (\ref{af1}), (\ref{as10}) and (\ref{as11}), we arrive at
the desired result (\ref{c3}). \endpf

\ms

Finally, similar to \cite[Lemma 3.3]{Zh2}, we have the following
result.

\bl\label{l3}
 Let $a^{jk}\in C^1(\cl{\O})$ satisfy (\ref{s3}), and
 $g\=(g^1, \cds, g^n): \dbR_t\t\dbR_x^n\rightarrow \dbR^n$ be
 a vector field of class $C^1$. Then for any $w\in C^2(\dbR_t\t\dbR_x^n;\;\dbC)$, we
 have
 \bel{df}
 \3n\3n\3n\ba{ll}
 &\displaystyle
 -\sum_{k=1}^n\left[(g\cdot\overline{\n w})\sum_{j=1}^n
 a^{jk}w_j+(g\cdot{\n w})\sum_{j=1}^n
 a^{jk}\overline w_j-g^k\left(|w_s|^2+\sum_{i, l=1}^na^{jl}w_j\overline w_l\right)\right]_k\\
 \noalign{\ms}
 &\displaystyle
 =-\left[w_{ss}+\sum_{j,k=1}^n(a^{jk}w_j)_k\right]g\cdot\overline{\n w}-\overline{\left(w_{ss}+\sum_{j,k=1}^n(a^{jk}w_j)_k\right)} g\cdot \n w\\
 \ns&\ds\q+(w_sg\cdot\overline{\n w}+\overline w_{s}g\cdot\n
w)_s-(w_sg_s\cdot\overline {\n w}+\overline w_sg\cdot\n w)\\
\ns&\ds\q +(\n\cdot
 g)|w_s|^2-2\sum_{j,k,l=1}^na^{jk}w_j\overline w_l\frac{\pa
 g^l}{\pa x_k}+\sum_{j,k=1}^nw_j\overline w_k\n\cd( a^{jk}g).
 \ea
 \ee
\el

{\it Proof.} On the one hand, we have
 \bel{df1}\ba{ll}\ds
w_{ss}g\cdot\overline{\n w}+\overline w_{ss}g\cdot\n w\\
\ns\ds=(w_sg\cdot\overline{\n w}+\overline w_{s}g\cdot\n
w)_s-(w_sg_s\cdot\overline {\n w}+\overline w_sg\cdot\n w)\\
\ns\ds\q-\sum_{j=1}^n(g^j|w_s|^2)_j+(\n\cdot g)|w_s|^2.
 \ea\ee
On the other hand, by (\ref{s3}), we have
 \bel{df2}\ba{ll}\ds
\sum_{j,k=1}^n(a^{jk}w_j)_kg\cdot\overline{\n w}+\overline{\sum_{j,k=1}^n(a^{jk}w_j)_k} g\cdot \n w\\
\ns\ds=\sum_{j,k=1}^n\[a^{jk}w_jg\cdot\overline{\n
w}+a^{jk}\overline w_j g\cdot \n
w\]_k-2\sum_{j,k,l=1}^na^{jk}w_j\overline w_l\frac{\pa
 g^l}{\pa x_k}\\
 \ns\ds\q-\sum_{j,k,l=1}^n(a^{jk}g^lw_j\overline w_k)_l+\sum_{j,k=1}^nw_jw_k\n\cd(
 a^{jk}g).
 \ea\ee
Combining (\ref{df1})--(\ref{df2}), we get the desired result.\endpf

\ms

\section{Global Carleman estimate for elliptic equations \\
without inhomogeneous boundary condition}\label{ss5}

In this section, we shall derive a global Carleman estimate for
elliptic equations with nonhomogeneous and complex Neumann-like
boundary condition.

Denote
$$
X=(-2,2)\t \O,\q \Si=(-2,2)\t\pa\O,\q Y=(-1,1)\t\O,\q
Z=(-2,2)\t\G_0.
$$
Let us consider the following  elliptic equation:
 \bel{s1}
 \left\{\ba{ll}
 \displaystyle
  z_{ss}+\sum_{j,k=1}^n\big(a^{jk}z_j\big)_k=z^0 & \mbox{ in } (-2,2)\t\O, \\
  \ns
 \displaystyle
\sum_{j,k=1}^na^{jk}z_j\nu_k=0& \mbox{ on }(-2,2)\t \pa\O\setminus\G_0,\\
  \ns
 \displaystyle
\sum_{j,k=1}^na^{jk}z_j\nu_k-ia(x)z_s=a(x)z^1& \mbox{ on
}(-2,2)\t\G_0.
 \ea\right.
 \ee
We now show the following Carleman estimate.

 \bt\label{t1} Under the assumptions in Theorem \ref{0t1}, there exists a constant $C>0$ such that, for any
$\e>0$, any
 solution $z\in C((-2,2);\;H^1(\O))\cap C^1((-2,2);\;L^2(\O))$ of system (\ref{s1})
 satisfies
  \bel{a2}\ba{ll}\ds
 ||z||_{H^1(Y)}\le
 Ce^{C\e}\[||z^0||_{L^2(X)}+||z^1||_{L^2(\Si)}+||z||_{L^2(Z)}+||z_s||_{L^2(Z)}\]\\
 \ns\ds\qq\qq\qq+Ce^{-2/\e}||z||_{H^1(X)}.
  \ea\ee
 \et

\br\label{remk1} For the general case of $t\in (T_1,T_2)$ with
$T_1,T_2\in\dbR$. By setting
 $s=t-{{T_2+T_1}\over 2}$, one deduces that
 $$s\in (-\a,\a), \q \a\= {{T_2-T_1}\over2}.$$
Then by scaling, one need consider only the case of (\ref{s1}).
 \er

\ms

{\it Proof.} We divide the proof into several steps.
\ms

{\it Step 1.} Note that there is no boundary condition for $z$ at
$s=\pm2$. Therefore, we need to introduce a cut-off function
$\varphi=\varphi(s)\in C_0^\i(-b,b)\subset C_0^\i(\dbR)$ such that
 \bel{as7}\left\{\ba{ll}\ds
 0\le\varphi(s)\le1 \q &|s|<b,\\
 \ns\ds \varphi(s)=1,&|s|\le b_0,
  \ea\right.\ee
where $b_0$ and $b$ (satisfying $1<b_0<b<2$) will be given later.
Put
 \bel{as8}
\hat z=\varphi z.
 \ee
Then, noting that $\varphi$ does not depend on $x$, by (\ref{s1}),
it follows
 \bel{as9}\left\{\ba{ll}\ds
\hat z_{ss}+\sum_{j,k=1}^n\big(a^{jk}\hat z_j\big)_k=\varphi_{ss}z+2\varphi_sz_s+\varphi z^0& \mbox{ in } (-2,2)\t\O, \\
 \ns
 \ds
\sum_{j,k=1}^na^{jk}\hat z_j\nu_k=0& \mbox{ on } (-2,2)\t \pa\O\setminus\G_0,\\
  \ns
 \ds
\sum_{j,k=1}^na^{jk}\hat z_j\nu_k-ia(x)\hat z_s=-ia(x)\varphi_sz+
a(x)\varphi z^1& \mbox{ on } (-2,2)\t\G_0.
 \ea\right.
 \ee

\ms

{\it Step 2.} Put
 \bel{as2}
b\=\sqrt{1+{1\over\mu}\ln \left[(2+e^\mu)e^{{\mu\hat\psi(x)\over
||\hat\psi||_{L^\i(\O)}}}\right]},\qq b_0\=\sqrt{b^2-{1\over\mu}\ln
\left[({1+e^\mu})e^{{\mu\hat\psi(x)\over
||\hat\psi||_{L^\i(\O)}}}\right]},
 \ee
where $\mu>\ln 2$, $\hat\psi(x)\in C^2(\oO)$ is given by Lemma
\ref{0l4}. It is easy to see that
 \bel{as3}
1<b_0<b\le2.
 \ee
Put
 \bel{as4}
 \psi=\psi(s,x)\=-{\hat\psi(x)\over ||\hat\psi||_{L^\i(\O)}}+b^2-s^2.
 \ee
It is easy to check that
 \bel{as5}\left\{\ba{ll}\ds
 \phi(s,\cd)\ge 2+e^\mu , &\hb{ for any $s$ satisfying }|s|\le 1,\\
 \ns\ds \phi(s,\cd)\le 1+e^{\mu}, &\hb{ for any $s$ satisfying }b_0\le |s|\le
 b.
 \ea\right.\ee
On the other hand, by (\ref{as4}) and Lemma \ref{0l4}, we find
 \bel{pb}
h\=|\n\psi|={1\over ||\hat\psi||_{L^\i(\O)}}|\n\hat\psi(x)|>0,\qq
\hb{in }\cl\O.
 \ee

Next, recalling that $a^{jk}$ satisfying (\ref{s4}) and by
(\ref{pb}), we conclude that there exists a $\mu_0>1$, for any
$\mu\ge\mu_0$, there exists $\l_0(\mu)>1$ such that for any $\l\ge
\l_1$, it holds
 \bel{pb8}\ba{ll}\ds
\hb {The right hand side of }
(\ref{c3})\\
\ns\ds\ge\l\mu^2\phi\sum_{j,k=1}^na^{jk}\psi_j\psi_k\(|v_s|^2+\sum_{j,k=1}^na^{jk}v_j\ov_k\)+\l^3\mu^4\phi^3\|\sum_{j,k=1}^na^{jk}\psi_j\psi_k\|^2|v|^2\\
\ns\ds\ge\l\mu^2\b
h^2\phi\(|v_s|^2+\sum_{j,k=1}^na^{jk}v_j\ov_k\)+\l^3\mu^4\b^2h^4\phi^3|v|^2.
 \ea\ee
Now, integrating inequality (\ref{c3}) (with $u$ replaced by $\hat
z$) in $(-b,b)\t \O$, recalling that $\varphi$ vanishes near $s=\pm
b$, and by (\ref{as9}) and (\ref{pb8}), one arrives at
 \bel{as12}\ba{ll}\ds
\l\mu^2\int_{-b}^b\int_{\O}\phi(|\n
v|^2+|v_s|^2)dxds+\l^3\mu^4\int_{-b}^b\int_{\O}\phi^3|v|^2dxds\\
\ns\ds\le C\Big\{\int_{-b}^b\int_{\O}\th^2|
\varphi_{ss}z+2\varphi_sz_s+\varphi
z^0|^2dxds+\int_{-b}^b\int_{\pa\O} V\cd\nu dxds\Big\}.
 \ea\ee
Recalling that $v=\th \hat z$, by (\ref{as6}), we get
 \bel{as14}
{1\over C}\th^2(|\n \hat z|^2+\l^2\mu^2\phi^2|\hat z|^2)\le |\n
v|^2+\l^2\mu^2\phi^2|v|^2\le C\th^2(|\n \hat z|^2+\l^2\mu^2\phi^2|
\hat z|^2).
 \ee
Therefore, by (\ref{as12}) and (\ref{as14}), we end up with
 \bel{as15}\ba{ll}\ds
\l\mu^2\int_{-b}^b\int_{\O}\th^2\phi(|\n
\hat z|^2+|\hat z_s|^2)dxds+\l^3\mu^4\int_{-b}^b\int_{\O}\th^2\phi^3|\hat z|^2dxds\\
\ns\ds\le C\Big\{\int_{-b}^b\int_{\O}\th^2|
\varphi_{ss}z+2\varphi_sz_s+\varphi
z^0|^2dxds+\int_{-b}^b\int_{\pa\O} V\cd\nu dxds\Big\}.
 \ea\ee

\ms

{\it Step 3.} We now estimate $\ds\int_{-b}^b\int_{\pa\O} V\cd\nu
dxds$. By (\ref{asf}) and nothing that $v=\th\hat z$, it follows
  \bel{as17}\3n\ba{ll}\ds
\int_{-b}^b\int_{\pa\O} V\cd\nu dxds=\sum_{k=1}^n\int_{-b}^b\int_{\pa\O} V_k\nu_k dxds\\
\ns\ds=\sum_{j,k,j',k'=1}^n\int_{-b}^b\int_{\pa\O} \Big\{-2a^{jk}\ell_j\nu_k|v_s|^2+2\ell_sa^{jk}\nu_k(\overline v_jv_s+v_j\overline v_s)-\Psi a^{jk}\nu_k(v_j\ov+\ov_jv)\\
 \ns\ds\q
+\(2a^{jk'}a^{j'k}-a^{jk}a^{j'k'}\)\ell_j(v_{j'}\ov_{k'}+\ov_{j'}v_{k'})\nu_k+a^{jk}\nu_k(2A\ell_j+\Psi_j-2a\ell_j\ell_t)|v|^2\Big\}dxds.
 \ea\ee

Note that, by (\ref{as}) and (\ref{as6}), we know that
 \bel{k1}
\sum_{j,k=1}^na^{jk}\ell_j\nu_k=\l\mu\phi
\sum_{j,k=1}^na^{jk}\psi_j\nu_k=-{\l\mu\phi\over
||\hat\psi||_{L^\i(\O)}} \sum_{j,k=1}^na^{jk}\hat\psi_j\nu_k\ge 0,\q
\hb{ on }\pa\O\setminus\G_0.
 \ee
Hence, recalling that $v=\th\hat z$, we have
 \bel{v1}\ba{ll}\ds
-\sum_{j,k=1}^n\int_{-b}^b\int_{\pa\O}a^{jk}\ell_j\nu_k|v_s|^2dxds\le
C\l\mu\int_{-b}^b\int_{\G_0}\phi|v_s|^2dxdt\\
\ns\ds\le Ce^{C\l}\int_{-b}^b\int_{\G_0}(|\hat z_s|^2+|\hat
z|^2)dxds.
 \ea\ee
Next, using $v=\th\hat z$ again, noting that $\hat z$ vanishes near
$s=\pm b$, by (\ref{a1}) and (\ref{as9}), we have
 \bel{v2}\ba{ll}\ds
\sum_{j,k=1}^n\int_{-b}^b\int_{\pa\O}\ell_sa^{jk}\nu_k(\overline v_jv_s+v_j\overline v_s)dxds-\sum_{j,k=1}^n\int_{-b}^b\int_{\pa\O}\Psi a^{jk}\nu_k(v_j\ov+\ov_jv)dxds \\
\ns\ds=\sum_{j,k=1}^n\int_{-b}^b\int_{\pa\O}\th^2\ell_sa^{jk}\nu_k(\overline
{\hat z_j}\hat z_s+\hat z_j\overline {\hat
z_s})dxds\\
\ns\ds\q+\sum_{j,k=1}^n\int_{-b}^b\int_{\pa\O}\th^2(\ell_s^2-\Psi)a^{jk}\nu_k(\overline
{\hat z_j}\hat z+\hat z_j\overline {\hat
z})dxds\\
\ns\ds\q+\sum_{j,k=1}^n\int_{-b}^b\int_{\pa\O}\th^2\ell_sa^{jk}\ell_j\nu_k(\hat
z\overline{\hat z_s}+\overline {\hat z}\hat
z_s)dxds\\
\ns\ds\q+2\sum_{j,k=1}^n\int_{-b}^b\int_{\pa\O}\th^2(\ell_s^2-\Psi)a^{jk}\ell_j\nu_k|\hat
z|^2dxds\\
\ns\ds=\int_{-b}^b\int_{\G_0}a(x)\th^2\ell_s\[i\varphi_s( \hat
z_s\overline {z}-\overline{\hat z_s}z)+\varphi(\overline
{\hat z_s}z^1+\hat z_s\overline {z^1})\]dxds\\
\ns\ds\q+\int_{-b}^b\int_{\G_0}a(x)\th^2(\ell_s^2-\Psi)\[i(\hat z_s\overline{\hat z}-\overline{\hat z_s}\hat z)+\varphi(\overline{z^1}\hat z+z^1\overline{\hat z})\]dxds\\
\ns\ds\q+\sum_{j,k=1}^n\int_{-b}^b\int_{\pa\O}(\th^2\ell_sa^{jk}\ell_j\nu_k|\hat
z|^2)_sdxds-\sum_{j,k=1}^n\int_{-b}^b\int_{\pa\O}\th^2(\ell_{ss}+2\Psi)a^{jk}\ell_j\nu_k|\hat
z|^2dxds\\
 \ns\ds\le Ce^{C\l}\[\int_{-b}^b\int_{\G_0}(|\hat
z_s|^2+|\varphi_s z|^2+|\varphi
z^1|^2)dxds+\int_{-b}^b\int_{\pa\O}|\hat z|^2dxds\].
 \ea\ee
Further, by (\ref{k1}), and noting that $v=\th \hat z $, we get
 \bel{as213}\ba{ll}\ds
\sum_{j,k,j',k'=1}^n\int_{-b}^b\int_{\pa\O}\(2a^{jk'}a^{j'k}-a^{jk}a^{j'k'}\)\ell_j(v_{j'}\ov_{k'}+\ov_{j'}v_{k'})\nu_k
dxds\\
\ns\ds=\sum_{j,k,j',k'=1}^n\int_{-b}^b\int_{\pa\O}a^{jk}\ell_j\nu_ka^{j'k'}(v_{j'}\ov_{k'}+\ov_{j'}v_{k'})dxds\\
\ns\ds\le Ce^{C\l}\sum_{j',k'=1}^n
\int_{-b}^b\int_{\pa\O\setminus\G_0}a^{j'k'}(v_{j'}\ov_{k'}+\ov_{j'}v_{k'})dxds\\
\ns\ds\le Ce^{C\l}
\int_{-b}^b\int_{\pa\O\setminus\G_0}\sum_{j,k=1}^n\[a^{jk}\hat
z_j\overline {\hat z_k}+a^{jk}\ell_j\ell_k|\hat z|^2\]dxds.
 \ea\ee
Combining (\ref{as17}), (\ref{v1})-- (\ref{as213}),  we obtain
 \bel{c5}\ba{ll}\ds
\int_{-b}^b\int_{\pa\O} V\cd\nu dxds\le
Ce^{C\l}\[\int_{-b}^b\int_{\G_0}(|\hat z_s|^2+|\varphi_sz|^2+|\varphi z^1|^2)dxds\\
\ns\ds\qq\qq\qq\qq\qq\qq+
\int_{-b}^b\int_{\pa\O\setminus\G_0}|\n\hat
z|^2dxds+\int_{-b}^b\int_{\pa\O}|\hat z|^2dxds\].
 \ea\ee

\ms

{\it Step 4.} Let us estimate $\ds \int_{-b}^b\int_{\pa\O}|\hat
z|^2dxds$ and $\ds \int_{-b}^b\int_{\pa\O\setminus\G_0}|\n\hat
z|^2dxds$.

Firstly, by trace theory and Poinc\'are inequality, noting that
$\hat z$ vanishes near $s=\pm b$, we have
 \bel{2pe2}\ba{ll}\ds
\int_{-b}^b\int_{\pa\O} |\hat z|^2dxds\le
C\int_{-b}^b\int_{\O}(|\hat z|^2+|\n \hat z|^2) dxds\le
C\int_{-b}^b\int_{\O}(|\hat z_s|^2+|\n \hat z|^2) dxds.
 \ea\ee

Next, we choose a $g\in C^1(\oO;\dbR)$ such that $g=\nu$ on $\pa\O$.
Integrating (\ref{df}) (in Lemma \ref{l3}) in $(-b,b)\t\O$, with $w$
replaced by $\hat z$, using integrating by parts, and noting $\hat
z(-b)=\hat z(b)=0$, by (\ref{as9}) and using Poinc\'are inequality,
we have
 \bel{df5}\ba{ll}\ds
 -\sum_{k=1}^n\int_{-b}^b\int_{\pa\O}\[(g\cdot\overline{\n
 \hat z})\sum_{j=1}^n
 a^{jk}\hat z_j\nu_k+(g\cdot{\n \hat z})\sum_{j=1}^n
 a^{jk}\overline {\hat z}_j\nu_k\]dxds\\
 \ns\ds\q+\int_{-b}^b\int_{\pa\O}\(|\hat z_s|^2+\sum_{j, l=1}^na^{jl}\hat z_j\overline {\hat z}_l\)dxds\\
 \noalign{\ss}
 \displaystyle
 =-\int_{-b}^b\int_\O\[\(\hat z_{ss}+\sum_{j,k=1}^n(a^{jk}\hat z_j)_k\)g\cdot\overline{\n \hat z}+\overline{\(\hat z_{ss}+\sum_{j,k=1}^n(a^{jk}\hat z_j)_k\)} g\cdot \n \hat z\]dxds\\
 \ns\ds\q-\int_{-b}^b\int_\O(\hat z_sg_s\cdot\overline {\n \hat z}+\overline {\hat z}_sg\cdot\n \hat z)dxds\\
\ns\ds\q +\int_{-b}^b\int_\O\[(\n\cdot
 g)|\hat z_s|^2-2\sum_{j,k,l=1}^na^{jk}\hat z_j\overline {\hat z}_l\frac{\pa
 g^l}{\pa x_k}+\sum_{j,k=1}^n\hat z_j\overline {\hat z}_k\n\cd( a^{jk}g)\]dxds\\
 \ns\ds\le C\int_{-b}^b\int_{\O}\[|\varphi_{ss}z+2\varphi_sz_s+\varphi z^0|^2+|\hat z_s|^2+|\n \hat
z|^2\]dxds\\
\ns\ds\le C\int_{-b}^b\int_{\O}(|z^0|^2+|\hat z_s|^2+|\n \hat
z|^2)dxds.
 \ea\ee
By (\ref{s4}), (\ref{as9}) and (\ref{df5}), we have
 \bel{dff5}\ba{ll}\ds
\int_{-b}^b\int_{\pa\O}(|\hat z_s|^2+\b|\n\hat
z|^2)dxds\le\int_{-b}^b\int_{\pa\O}\(|\hat z_s|^2+\sum_{j,
l=1}^na^{jl}\hat z_j\overline {\hat z}_l\)dxds\\
\ns\ds\le C\int_{-b}^b\int_{\O}(|\hat z_s|^2+|\n \hat
z|^2+|z^0|^2) dxds\\
\ns\ds\q+\int_{-b}^b\int_{\G_0}a(x)\[(g\cdot\overline{\n
 \hat z})(i\hat z_s-i\varphi_s z+\varphi z^1)+(g\cdot{\n \hat z})(-i\overline{\hat z_s}+i\varphi_s\overline
 z+\varphi\overline{z^1})\]dxds\\
 \ns\ds\le C\int_{-b}^b\int_{\O}(|\hat z_s|^2+|\n \hat
z|^2+|z^0|^2) dxds\\
\ns\ds\q+\d\int_{-b}^b\int_{\G_0}|\n\hat
z|^2dxds+C(\d)\int_{-b}^b\int_{\G_0}(|\varphi_s z|^2+|\hat
z_s|^2+|\varphi z^1|^2)dxds
 \ea\ee
where $0<\d<\b$ is small, then
 \bel{df7}\ba{ll}\ds
\int_{-b}^b\int_{\pa\O\setminus\G_0}|\n\hat z|^2dxds\\
\ns\ds\le C\[\int_{-b}^b\int_{\O}(|\hat z_s|^2+|\n \hat
z|^2+|z^0|^2) dxds+\int_{-b}^b\int_{\G_0}(|\hat z|^2+|\hat
z_s|^2+|z^1|^2)dxds\].
 \ea\ee

Finally, by multiplying $\overline {\hat z}$ and $\hat z$ on the
first equation of (\ref{as9}), respectively, using integrating by
parts, by (\ref{s4}) and using Poinc\'are inequality, we get
 \bel{2pe3}\ba{ll}\ds
 2\int_{-b}^b\int_{\O}(|\hat z_s|^2+\b |\n\hat z|^2)dxds\\
 \ns\ds
\le\int_{-b}^b\int_{\O}\(2|\hat z_s|^2+\sum_{j,k=1}^na^{jk}(\hat z_j\overline{\hat z_k}+\overline{\hat z_j}\hat z_k)\)dxds\\
\ns\ds=\sum_{j,k=1}^n\int_{-b}^b\int_{\pa\O}\(\overline {\hat
z}a^{jk}\hat z_j\nu_k+\hat za^{jk}\overline{\hat z_j}\nu_k\)
dxds-\int_{-b}^b\int_{\O}\varphi(z^0\overline{\hat z}+\overline z^0 \hat z)dxds\\
\ns\ds\le C\[\int_{-b}^b\int_{\G_0}(|\hat z|^2+|\hat
z_s|^2+|z^1|^2)dxds\]+{1\over\e^*}\int_{-b}^b\int_{\O}
|z^0|^2dxds+\e^*\int_{-b}^b\int_{\O} |\hat z|^2dxds\\
\ns\ds\le C\[\int_{-b}^b\int_{\G_0}(|\hat z|^2+|\hat
z_s|^2+|z^1|^2)dxds\]\\
\ns\ds\q+{1\over\e^*}\int_{-b}^b\int_{\O}
|z^0|^2dxds+C\e^*\int_{-b}^b\int_{\O}|\hat z_s|^2dxds.
  \ea\ee
Taking $\e^*={1\over C}$ small enough, and combining (\ref{2pe2}),
(\ref{df7}) and (\ref{2pe3}), we get
 \bel{pe1}\ba{ll}\ds
 \int_{-b}^b\int_{\pa\O} |\hat z|^2dxds+\int_{-b}^b\int_{\pa\O\setminus\G_0}|\n\hat z|^2dxds\\
 \ns\ds\le C\[\int_{-b}^b\int_{\G_0}(|\hat z|^2+|\hat
z_s|^2+|z^1|^2)dxds+\int_{-b}^b\int_{\O} |z^0|^2dxds\].
 \ea\ee

By (\ref{c5}) and (\ref{pe1}), and noting that $\hat z=\varphi z$,
we obtain
 \bel{as2a5}\ba{ll}\ds
\int_{-b}^b\int_{\pa\O} V\cd\nu dxds\\
\ns\ds\le
Ce^{C\l}\[\int_{-b}^b\int_{\O}|z^0|^2dxds+\int_{-b}^b\int_{\G_0}
(|z|^2+|z_s|^2+|z^1|^2) dxds\].
 \ea\ee

\ms

{\it Step 5.} Combing (\ref{as15}), (\ref{df7}) and (\ref{as2a5}),
we end up with
  \bel{as216}\ba{ll}\ds
\l\mu^2\int_{-b}^b\int_{\O}\th^2\phi(|\n
 z|^2+|z_s|^2)dxds+\l^3\mu^4\int_{-b}^b\int_{\O}\th^2\phi^3|z|^2dxds\\
\ns\ds\le C\int_{-b}^b\int_{\O}\th^2\phi^3|
 \varphi_{ss}z+2\varphi_sz_s+\varphi z^0|^2dxds\\
 \ns\ds\q+Ce^{C\l}\[\int_{-b}^b\int_{\O}|z^0|^2dxdt+\int_{-b}^b\int_{\G_0}
(|z|^2+|z_s|^2+|z^1|^2) dxds\].
 \ea\ee
Denote $c_0=2+e^{\mu}>1$, and recall (\ref{as2}) for $b_0\in (1,b)$.
Fixing the parameter $\mu$ in (\ref{as216}), using (\ref{as7}) and
(\ref{as5}), one finds
 \bel{as217}\ba{ll}\ds
\l e^{2\l c_0}\int_{-1}^1\int_{\O}(|\n
 z|^2+|z_s|^2+|z|^2)dxds\\
\ns\ds\le
Ce^{C\l}\Big\{\int_{-2}^2\int_{\O}|z^0|^2dxds+\int_{-2}^2\int_{\pa
\O}|z^1|^2 dxds+\int_{-2}^2\int_{\G_0}(|z|^2+|z_s|^2)
dxds\Big\}\\
\ns\ds\q+Ce^{2\l(c_0-1)}\int_{(-b,-b_0)\bigcup (b_0,
b)}\int_{\O}(|z|^2+|z_s|^2)dxds.
 \ea\ee
From (\ref{as217}), one concludes that there exists an $\e_2>0$ such
that the desired inequality (\ref{a2}) holds for $\e\in (0,\e_2]$,
which, in turn, implies that it holds for any $\e>0$. This completes
the proof of Theorem \ref{t1}.\endpf

\ms

\section{Proof of Theorem \ref{0t2}}\label{ss6}

In this section, we will prove the existence and the estimate of the
norm of the resolvent $(A-\l I)^{-1}$ when $
 \Re\l\in\[-e^{-C|\Im\l|}/ C,0\]$.

\ms

{\it Proof. } We divide the proof into two steps.

\ms

{\it Step 1. } First, let $f=(f^0,f^1)\in H$, and $u=(u^0,u^1)\in
D(A)$ with the boundary condition
$\ds\sum_{j,k=1}^na^{jk}u^0_{j}\nu_k\|_{ \pa\O\setminus\G_0}=0 ,\
\(\sum_{j,k=1}^na^{jk}u^0_{j}\nu_k+au^1\)\|_{\G_0}=0 $.

Then, the following equation
 \bel{6a1}
(A-\l I)u=f
 \ee
is equivalent to
 \bel{6a2}\left\{\ba{ll}\ds
 -\l u^0+u^1=f^0,\\
 \ns\ds
 \sum_{j,k=1}^n(a^{jk}u^0_{j})_k-\l u^1=f^1.\ea\right.\ee

Hence, by substituting $u^1$ by $u^0$ in the second equation of
(\ref{6a2}) and with the boundary condition, we have
 \bel{6a3}\left\{\ba{ll}\ds
\sum_{j,k=1}^n(a^{jk}u^0_{j})_k-\l^2 u^0=\l f^0+f^1 & \hb { in } \O,\\
\ns\ds \sum_{j,k=1}^na^{jk}u^0_{j}\nu_k=0 &\hb{ on } \pa\O\setminus\G_0,\\
\ns\ds \sum_{j,k=1}^na^{jk}u^0_{j}\nu_k+a\l u^0=-af^0 &\hb{ on } \G_0,\\
\ns\ds u^1=f^0+\l u^0 & \hb { in } \O.
 \ea\right.\ee
Put
 \bel{v}
v=e^{i\l s}u^0.
 \ee
It is easy check that $v$ satisfying the following equation:
 \bel{6a4}\left\{\ba{ll}\ds
v_{ss}+\sum_{j,k=1}^n(a^{jk}v_j)_k=(\l f^0+f^1)e^{i\l s} & \hb { in } \dbR\t \O,\\
\ns\ds \sum_{j,k=1}^na^{jk}v_j\nu_k=0& \hb { on } \dbR\t \pa\O\setminus\G_0,\\
\ns\ds \sum_{j,k=1}^na^{jk}v_j\nu_k-ia v_s=-af^0e^{i\l s}& \hb { on
}\dbR\t \G_0.
 \ea\right.\ee

\ms

{\it Step 2. } By (\ref{v}) and Remark \ref{remk1}, we have the
following estimates.
 \bel{6a7}\left\{\ba{ll}\ds
|u^0|_{H^1(\O)}\le Ce^{C|\Im\l|}|v|_{H^1(Y)},\\
\ns\ds |v|_{H^1(X)}\le C(|\l|+1)e^{C|\Im\l|}|u^0|_{H^1(\O)},\\
 \ns\ds |v|_{L^2(Z)}\le Ce^{C|\Im\l|}|u^0|_{L^2(\G_0)},\q
|v_s|_{L^2(Z)}\le C|\l|e^{C|\Im\l|}|u^0|_{L^2(\G_0)}.
 \ea\right.\ee

Now, we apply $v$ to Theorem \ref{t1}, and combining (\ref{6a7}), we
have
 \bel{6a5}
|u^0|_{H^1(\O)}\le
Ce^{C|\Im\l|}\[|f^0|_{H^1(\O)}+|f^1|_{L^2(\O)}+|u^0|_{L^2(\G_0)}\].
 \ee
On the other hand, we multiplier (\ref{6a2}) by $\overline u^0$,
integrate it on $\O$, we get
 \bel{6a8}\ba{ll}\ds
\int_{\O}\(-\sum_{j,k=1}^n(a^{jk}u^0_{j})_k+\l^2u^0\)\cd\overline u^0dx\\
\ns\ds=\l^2|u^0|^2_{L^2(\O)}+\sum_{j,k=1}^n\int_{\O}
a^{jk}u^0_{j}\overline u^0_{k}dx-\sum_{j,k=1}^n\int_{\pa\O}
a^{jk}u^0_{j}\nu_k\overline u^0dx\\
\ns\ds=\l^2|u^0|^2_{L^2(\O)}+\sum_{j,k=1}^n\int_\O
a^{jk}u^0_{j}\overline u^0_{k}dx+\int_{\pa\O}(a\l u^0+af^0)\overline
u^0dx.
 \ea\ee
By taking the imaginary part, we find,
 \bel{6a9}\ba{ll}\ds
|\Im\l|\int_{\pa\O} a|u^0|^2dx\\
\ns\ds\le
\|-\sum_{j,k=1}^n(a^{jk}u^0_{j})_k+\l^2u^0\|_{L^2(\O)}|u^0|_{L^2(\O)}\\
\ns\ds\qq+2|\Im\l||\Re\l||u^0|^2_{L^2(\O)}+C|f^0|_{L^2(\pa\O)}|\sqrt{a}u^0|_{L^2(\pa\O)}\\
\ns\ds\le C\[|(\l
f^0+f^1)|_{L^2(\O)}|u^0|_{L^2(\O)}+|\Im\l||\Re\l||u^0|^2_{L^2(\O)}+|f^0|_{H^1(\O)}|u^0|_{H^1(\O)}\]
 \ea\ee
Hence, combining (\ref{6a5}) and (\ref{6a9}), we have
 \bel{6a10}\ba{ll}\ds
|u^0|_{H^1(\O)}\le
Ce^{C|\Im\l|}\[|f^0|_{H^1(\O)}+|f^1|_{L^2(\O)}+|\Im\l||\Re\l||u^0|_{H^1(\O)}\].
 \ea\ee
Therefore, we take
 $$
Ce^{C|\Im\l|}|\Im\l||\Re|\le{1\over2},
$$
which holds, as soon as $\ds |\Re\l|\le -e^{C_0|\Im\l|}/C_0$ for
some $C_0>0$. Then, we have
  \bel{6a11}\ba{ll}\ds
|u^0|_{H^1(\O)}\le Ce^{C|\Im\l|}(|f^0|_{H^1(\O)}+|f^1|_{L^2(\O)}).
 \ea\ee
Recalling that $u^1=f^0+\l u^0$, we have
  \bel{6a12}\ba{ll}\ds
|u^1|_{L^2(\O)}\le |f^0|_{L^2(\O)}+|\l||u^0|_{L^2(\O)}\le
Ce^{C|\Im\l|}(|f^0|_{H^1(\O)}+|f^1|_{L^2(\O)}).
 \ea\ee
By (\ref{6a11})--(\ref{6a12}), we know that $A-\l I$ is injective.
Thus $A-\l I$ is bi-injective from $D(A)$ to $H$. And moreover,
 $$
||(A-\l I)^{-1}||_{\cL(H,H)}\le Ce^{C|\Im\l|},\q  \Re\l\in
(-e^{C|\Im\l|}/C,0),\qq |\l|\ge 1.
 $$
This completes the proof of Theorem \ref{0t2}.\endpf

\ms

\section{Proof of Theorem \ref{0t1}}\label{ss7}

In this section, we adapt the proof of \cite[Th\'eor\`eme 3]{B} (and
also the proof of \cite[Theorem 3]{C} on semigroups).

{\it Proof of Theorem \ref{0t1}.} By taking $\chi_1=\chi_2=I$,
$A=iB$ and $k=2$ in \cite[Th\'eor\`eme 3]{B}, we have
 \bel{w1}
\|\|e^{tA}u{1\over {(I-A)^2}}\|\|_H\le  \left({C\over \ln
(2+t)}\right)^2||u||_H,
 \ee
that is
 \bel{w2}
||e^{tA}u||_H\le   \left({C\over \ln (2+t)}\right)^2||u||_{D(A^2)}.
 \ee
By definition, $D(A)$ is the interpolate space  between $D(A^0)=H$
and $D(A^2)$. Since
 \bel{w3}
||e^{tA}u||_H\le C||u||_H.
 \ee
Then, combining (\ref{w2})--(\ref{w3}), by applying interpolation
theorem, we get the desired result.\endpf


\begin{thebibliography}{99}

\bibitem{B} N.~Burq,   \it D\'ecroissance de l'\'energie locale de l'\'equation des ondes pour le probl\'eme ext\'erieur et absence de r\'esonance au voisinagage du r\'eel,
\sl Acta Math., \rm  180 (1998), 1--29.

\bibitem{BH} N.~Burq and M.~Hitrik, \it Energy decay for damped wave equations on
partially rectangular domains, \sl Math. Res. Lett., \rm 14 (2007),
35--47.

\bibitem{BLR} C.~Bardos, G.~Lebeau and J.~Rauch,   \it Sharp sufficient conditions for the observation, control, and stabilization of waves from the boundary,
\sl SIAM I. Control Optim., \rm 30 (1992), 1024--1065.

\bibitem{C} H.~Christianson,   \it Applications of cutoff resolvent estimates to the wave equations,
\rm Preprint.

\bibitem{D} T.~Duyckaerts, \it Optimal decay rates of the energy of a
hyperbolic-parabolic system coupled by an interface, \sl Asymptot.
Anal., \rm 51 (2007), 17--45.

\bibitem{DZZ} T.~Duyckaerts, X.~Zhang and E.~Zuazua, \it On the
optimality of the observability inequalities for parabolic and
hyperbolic systems with potentials, \sl Ann. Inst. H. Poincar\'e
Anal. Non Lin\'eaire,  \rm 25 (2008), 1--41.

\bibitem{Fu0} X.~Fu,  \it A weighted identity for partial differential operators of second order
and its applications, \sl C. R. Math. Acad. Sci. Paris, \rm  342
(2006),  579--584.

\bibitem{Fu} X.~Fu,  \it Null controllability for the parabolic equations with a complex principal part, \rm Preprint.

\bibitem{FI} A.~V.~Fursikov and O.~Yu.~Imanuvilov, \sl
Controllability of Evolution Equations, \rm Lecture Notes Series 34,
Research Institute of Mathematics, Seoul National University, Seoul,
Korea, 1994.

\bibitem{LRS}  M.~M.~Lavrent'ev, V.~G.~Romanov, S.~P.~Shishat.skii,
\sl Ill-Posed Problems of Mathematics Physics and Analysis, \rm
Translated from the Russian by J.~R.~Schulenberger, Translations of
Mathematical Monographs, 64,  American Mathematical Society,
Providence, RI, 1986.

\bibitem{LR} G.~Lebeau  and L.~Robbiano,
\it Stabilisation de l'\'equation des ondes par le bord, \sl Duke
Math. J., \rm 86 (1997), 465--491.

\bibitem{LR1} Z.~Liu and B.~Rao,   \it Characterization of polynomial decay rate for the solution of linear evolution equation,
\sl Z. angew. Math. Phys., \rm  56 (2005),  630--644.

\bibitem{P1} K.-D.~Phung, \it Polynomial decay rate for the dissipative wave
equation, \sl J. Differential Equations, \rm 240 (2007),  92--124.

\bibitem{P2} K.-D.~Phung, \it Boundary stabilization for the wave equation in
a bounded cylindrical domain, \sl Discrete Contin. Dynam. Systems,
\rm 20 (2008), 1057--1093.

\bibitem{RZZ} J.~Rauch, X.~Zhang and E.~Zuazua, \it Polynomial decay
of a hyperbolic-parabolic coupled system, \sl J. Math. Pures Appl.,
\rm 84 (2005), no. 4, 407--470.

\bibitem{WW}G.~Wang and L.~Wang, \it The Carleman inequality and its
application to periodic optimal control governed by semilinear
parabolic differential equations, \sl J. Optim. Theory Appl., \rm
118 (2003),  249--461.

\bibitem{Zh2} X.~Zhang, \it Explicit observability estimate for the wave equation with
potential and its application,  \sl  R. Soc. Lond. Proc. Ser. A
Math. Phys. Eng. Sci., \rm {\bf 456} (2000), 1101--1115.

\bibitem{Z1} X.~Zhang, \it Explicit observability inequalities for
the wave equation with lower order terms by means of Carleman
inequalities, \sl SIAM J. Control Optim., \rm 39 (2000), 812--834.

\bibitem{ZZ2} X.~Zhang and E.~Zuazua, \it Long time behavior of a
coupled heat-wave system arising in fluid-structure interaction, \sl
Arch. Rat. Mech. Anal., \rm 184 (2007), 49--120.

\end{thebibliography}
\end{document}